\documentclass[12pt]{amsart}
\usepackage{graphicx}
\usepackage{amssymb}
\usepackage{amsfonts}
\usepackage{amsmath}
\usepackage{array}
\usepackage{rotating}
\usepackage{bm}
\usepackage[matrix,frame]{xypic}

\newtheorem{example}{Example}[section]

\newtheorem{theorem}[example]{Theorem}

\def\Proof{\noindent \it Proof -- \rm}
\def\qed{\hspace{3.5mm} \hfill \vbox{\hrule height 3pt depth 2 pt width 2mm}
\bigskip}

\def\FQSym{{\bf FQSym}}

\def\ssh{\Cup}

\def\sconc{\bullet}
\def\Std{{\rm Std}}

\def\<{\langle}
\def\>{\rangle}

\def\Q{\operatorname{\mathbb Q}}
\def\N{\operatorname{\mathbb N}}

\def\F{{\bf F}}
\def\S{{\bf S}}

\def\G{{\bf G}}

\def\SG{{\mathfrak S}}

\def\shuff#1#2{\mathbin{
\hbox{\vbox{ \hbox{\vrule \hskip#2 \vrule height#1 width 0pt
}%
\hrule}%
\vbox{ \hbox{\vrule \hskip#2 \vrule height#1 width 0pt
\vrule }%
\hrule}%
}}}

\def\shuf{{\mathchoice{\shuff{7pt}{3.5pt}}%
{\shuff{6pt}{3pt}}%
{\shuff{4pt}{2pt}}%
{\shuff{3pt}{1.5pt}}}}%
\def\shuffle{\,\shuf\,}

\def\fatter{\preceq}

\def\diam{{\rm diam}}

\def\carr{\blacktriangleright}

\title{Inversion of some series of free quasi-symmetric functions}

\author[F. Hivert, J.-C.~Novelli, and J.-Y.~Thibon]%
{Florent Hivert, Jean-Christophe Novelli, and Jean-Yves Thibon}

\address[Hivert]{LITIS, Universit\'e de Rouen ; Avenue de l'universit\'e ;
76801 Saint \'Etienne du Rouvray, France\\}

\address[Novelli and Thibon]{Universit\'e Paris-Est, Institut Gaspard Monge
\\ 5 Boulevard Descartes \\Champs-sur-Marne \\77454 Marne-la-Vall\'ee cedex 2 \\
France}

\email[Florent Hivert]{hivert@univ-rouen.fr}
\email[Jean-Christophe Novelli]{novelli@univ-mlv.fr}
\email[Jean-Yves Thibon]{jyt@univ-mlv.fr} 

\begin{document}

\begin{abstract}
We give a combinatorial formula for the inverses of the alternating sums
of free quasi-symmetric functions of the form $\F_{\omega(I)}$ where $I$
runs over compositions with parts in a prescribed set $C$. This
proves in particular three special
cases (no restriction, even parts, and all parts equal to 2) which were
conjectured by B. C. V. Ung in [Proc. FPSAC'98, Toronto].
\end{abstract}

\maketitle

\section{Introduction}

The algebra of Free Quasi-Symmetric Functions $\FQSym$ \cite{NCSF6}
is a graded algebra of noncommutative polynomials whose bases are parametrized
by permutations. Under commutative image, it is mapped onto Gessel's algebra
of quasi-symmetric functions, whence its name.

Quasi-symmetric functions generalize symmetric functions in a natural way,
and many classical results admit quasi-symmetric extensions or analogs.
However, very few results resembling symmetric series identities, like
those of Schur or Littlewood (see, e.g., \cite{Mcd}) are known.
In \cite{Ung}, B. C. V. Ung proves a quasi-symmetric analog of Schur's
identity, and conjectures three further combinatorial inversions of 
quasi-symmetric series, which are even stated at the level of $\FQSym$.

In this note, we prove a master identity, which consists in a
combinatorial formula for the inverses of the alternating sums
of free quasi-symmetric functions of the form $\F_{\omega(I)}$ where $I$
runs over compositions with parts in a prescribed set $C$.
Here $\F_\sigma$ denotes the standard basis of $\FQSym$
(mapped onto Gessel's fundamental basis), and $\omega(I)$ is the longest
permutation with descent composition $I$.
Ung's conjectures boil down to the following special cases : no restriction
on the parts, even parts, and all parts equal to 2.

{\footnotesize
{\it Acknowledgements.}
This project has been partially supported by the grant ANR-06-BLAN-0380.
The authors would also like to thank the contributors of the MuPAD project,
and especially of the combinat part, for providing the development environment
for their research (see~\cite{HTm} for an introduction to MuPAD-Combinat).
}

\section{Background and notations}

\subsection{Free quasi-symmetric functions}

Let $A$ be a totally ordered alphabet.
Recall that the standardized
$\Std(w)$ of a word $w\in A^*$ is
the permutation obtained by iteratively scanning $w$ from
left to right, and labelling $1,2,\ldots$ the occurrences of its
smallest letter, then numbering the occurrences of the next one, and
so on. Alternatively, $\sigma=\Std(w)^{-1}$ can be characterized as
the unique permutation of minimal length such that $w\sigma$ is a
nondecreasing word. For example, $\Std(bbacab)=341625$.

An elementary observation,
which is at the basis of the constructions of \cite{NCSF6}, is that
the noncommutative polynomials
\begin{equation}
\G_\sigma(A)=\sum_{w\in A^*;\Std(w)=\sigma}w
\end{equation}
span a subalgebra of $\Q\langle A\rangle$. 
When $A$ is
infinite, this subalgebra admits a natural Hopf algebra structure,
but this fact will not be needed here.
This is $\FQSym$, the algebra of {\em Free Quasi-Symmetric Functions}.

Let $\F_\sigma=\G_{\sigma^{-1}}$. 
The scalar product is defined
by
\begin{equation}
\<\F_\sigma\,,\,\G_\tau\>=\delta_{\sigma,\tau}\,.
\end{equation}
For a word $w$ on the alphabet $\{1,2,\ldots\}$, denote by $w[k]$ the word
obtained by replacing each letter $i$ by the integer $i+k$.
If $u$ and $v$ are two words, with $u$ of length $k$, one defines
the {\em shifted concatenation}
\begin{equation}
u\sconc v = u\cdot (v[k])
\end{equation}
and the {\em shifted shuffle}
\begin{equation}
u\ssh v= u\shuffle (v[k])\,.
\end{equation}
where $u\shuffle u'$ is the usual shuffle product on words.

The product formula in the $\F$ basis is
\begin{equation}
\F_\alpha\F_\beta=\sum_{\gamma\,\in\,\alpha\ssh\beta }\F_\gamma\,.
\end{equation}
The sum of the inverses of the permutations occuring in
$\alpha^{-1}\ssh\beta^{-1}$ is called
{\em convolution} and denoted by $\alpha \star \beta$.
Thus, $\FQSym$ provides a realization of the convolution
algebra of permutations studied in  \cite{Re,MR}.

\subsection{Descent classes and compositions}

Recall that the descent set of a permutation $\sigma$
is $D=\{i\,|\, \sigma(i)>\sigma(i+1)\}$.
If $\sigma\in\SG_n$ has descent set
$D = \{d_1<\ldots<d_k\}\subseteq [n-1]$, the {\it descent composition}
$I=C(\sigma)$ is the composition $I=(i_1,\ldots,i_{k+1})$ of $n$ defined by
$i_s=d_s-d_{s-1}$, where $d_0=0$ and $d_{k+1}=n$. 
The symbol $I\vDash n$ means that $I$ is a composition of $n$, and
$l(I)$ denotes the length of $I$.

The descent class $D_I=\{\sigma\in\SG_n\,|\, C(\sigma)=I\}$
has a unique element of minimal (resp. maximal) length
denoted by $\alpha(I)$ (resp. $\omega(I)$).
Actually, descent classes are intervals $D_I=[\alpha(I),\omega(I)]$ for the
left weak order on $\SG_n$ (see, {\it e.g.}, \cite{BB}).

The mirror image of a word $w=a_1a_2\cdots a_m$ 
is $\bar w=a_ma_{m-1}\cdots a_1$. We shall use this
notation for compositions and permutations as well.

Finally, the diameter of a descent class is the permutation
\begin{equation}\label{diam}
\diam(I) := \alpha(I) \omega(I)^{-1}
          = \alpha(I) \omega(\overline{I}).
\end{equation}

\subsection{A multiplicative basis of $\FQSym$}

The \emph{left-shifted concatenation} of words is
\begin{equation}
u\carr v = u[l]\cdot v \quad {\rm if}\ u\in A^k,\ v\in A^l\,,
\end{equation}
similar to the usual shifted concatenation
$\bullet$, but with the shift on the first factor.
The following basis is introduced in \cite{NCSF7}:
\begin{equation}
\label{eq:defS}
\S^\sigma := \sum_{\tau\leq\sigma}\G_\tau\,
\end{equation}
where $\leq$ is the left weak order.
It has the property
\begin{equation}
\label{eq:multS}
\S^\sigma =\S^{\sigma_1}\S^{\sigma_2}\cdots \S^{\sigma_r}
\end{equation}
whenever $\sigma=\sigma_1\carr\sigma_2\carr\cdots\carr\sigma_r$.

The  Moebius function of the left weak order is explicitely
known \cite{Ede,Bj,BB}, and gives in particular
\begin{equation}
\label{G2S}
\G_\sigma = \sum_{I\fatter C(\sigma^{-1})} (-1)^{l(I)-1}\S^{\alpha(I)\sigma}.
\end{equation}

\section{The main result}

\subsection{Ung's conjectures}

In~\cite{Ung}, Ung made the following conjectures.
The inverses of the series
\begin{eqnarray*}
  H_1 &=& \sum_I (-1)^{\ell(I)}\F_{\omega(I)}\\
  H_2 &=& \sum_{n\ge 0} (-1)^n\F_{\omega(2^n)}\\
  H_3 &=&  \sum_I (-1)^{\ell(I)}\F_{\omega(2I)}\\
\end{eqnarray*}
are as follows. For a permutation $\sigma$ of shape
$I$, let $\hat\sigma=\sigma\omega(I)^{-1}$. Then,
\begin{eqnarray*}
  H_1^{-1} &=& \sum_\alpha \G_{\hat\alpha}\\
  H_2^{-1} &=& \sum_\beta \G_{\hat\beta}\\
  H_3^{-1}&=&\sum_\gamma \G_{\hat\gamma}
\end{eqnarray*}
where $\alpha$ runs over all permutations, $\beta\in\SG_{2p}$ runs
over permutations of shape $2^{2p}$, and $\gamma\in\SG_{2p}$ runs over
permutations with descent set contained in $\{2,4,\ldots,2p-2\}$.

Taking into account (\ref{eq:defS}) and (\ref{diam}), 
we see that all three identities are of the form (\ref{general}) below,
with $E=\N^*$, $\{2\}$ and $2\N^*$, respectively.

\subsection{Generalization}

\begin{theorem}
\label{prodG}
Let $E$ be any subset of $\N^*$.
And let $C(E)$ be the set of all compositions with parts in this subset.
Then
\begin{equation}\label{general}
\left( \sum_{I\in C(E)} (-1)^{l(I)} \G_{\omega(I)} \right)^{-1} =
\sum_{K\in C(E)} \S^{\diam(K)}.
\end{equation}
\end{theorem}

\Proof
Thanks to~(\ref{G2S}), the statement to be proved is equivalent to
\begin{equation}
\left( \sum_{I\in C(E)} (-1)^{l(I)}
       \sum_{J\fatter\overline{I}} (-1)^{l(J)-1} \S^{\alpha(J)\omega(I)}
\right)
\left(\sum_{K\in C(E)} \S^{\alpha(K)\omega(\overline{K})}\right)
= 1,
\end{equation}
or, opening the parentheses,
\begin{equation}
\label{sommetot}
\sum_{I,K\in C(E)}
\sum_{J\fatter \overline{I}}
(-1)^{l(I)+l(J)-1} \S^{\alpha(J)\omega({I})}
\S^{\alpha(K)\omega(\overline{K})} = 1.
\end{equation}
Now,
\begin{equation}
\S^{\alpha(J)\omega({I})} \S^{\alpha(K)\omega(\overline{K})}
= \S^{\alpha(J') \omega(I')},
\end{equation}
where $I'=I\bullet \overline{K}$ and
      $J'=K \triangleright J$.
Note that $J'\fatter \overline{I'}$ and that
\begin{equation}
(-1)^{l(I)+l(J)-1} = - (-1)^{(l(I')+l(J')-1}.
\end{equation}
Now, given any non-empty permutation $\sigma$ obtained as a product
$\alpha(J)\omega(I)$ with $J\fatter\overline{I}$, it can be decomposed in
exactly two ways as a product
$\alpha(J)\omega({I})\carr\alpha(K)\omega(\overline{K})$:
either with $K=\emptyset$ or with $\alpha(K)\omega(\overline{K})$
corresponding to the last anticonnected permutation associated with the
decomposition of $\sigma$ into anticonnected permutations.
This comes from the fact that $\alpha(J)\omega(I)$ (with
$J\fatter\overline{I}$) is anticonnected iff $J=\overline{I}$.

Since the coefficients associated with these two decompositions are opposite,
such a permutation does not occur in the final result. Hence the result
reduces to the contribution of the empty permutation.
\qed

\section{Comments on Ung's other identities}

In \cite{Ung}, Ung proves quasi-symmetric analogs of Schur's identity
(for the sum of all Schur functions) and of Littlewood's identity (for
its inverse). In fact, these analogs may be formulated 
without further work  at the level
of $\FQSym$. 

The first identity is
\begin{equation}\label{qschur}
\sum_I F_I = \frac12\left[\prod_i
\frac{1+x_i}{1-x_i}-1\right]=\frac12[\lambda_1(X)\sigma_1(X)-1]
\end{equation}
where $\lambda_1$ (resp. $\sigma_1$) is the sum of the elementary
(resp. complete) symmetric functions. Interpreting the right-hand
side in the algebra of noncommutative symmetric functions,
we have
\begin{equation}\label{ncschur}
\frac12[\lambda_1(A)\sigma_1(A)-1]=
\frac12\left[\prod_i^\leftarrow(1+a_i)\prod_i^\rightarrow(1-a_i)^{-1}-1\right]
=\sum_{n\ge 0}H_n
\end{equation}
where 
\begin{equation}
H_n=\sum_{k=0}^{n-1}R_{1^k,n-k}\,.
\end{equation}
The commutative image of $R_{1^k,n-k}$ is the Schur function
$s_{n-k,1^k}$, whose quasi-symmetric expansion is easily found
to be
\begin{equation}
s_{n-k,1^k}=\sum_{I\vdash n,\ l(I)=k+1}F_I\,.
\end{equation}
But $R_{1^k,n-k}$ can also be interpreted as an element
of $\FQSym$,
\begin{equation}
R_{1^k,n-k}=\sum_{C(\sigma^{-1})=(1^k,n-k)}\F_{\sigma}
\end{equation}
so that (\ref{qschur}) means that each descent class contains
exactly one permutation whose inverse has a hook shape $(1^k,n-k)$.

The second identity is
\begin{equation}\label{qlit}
\left(\sum_I F_I\right)^{-1}=1+\sum_{I\vDash 2n+1}(-1)^{n+1}c_IF_I\,
\end{equation}
where $c_I$ is the number of permutations of shape $I$ whose inverse
has shape $(12^n)$. This formula is obtained by observing that
the inverse of $H=\sum_n H_n$ is the noncommutative hyperbolic tangent
of \cite{NCSF1}, that is
\begin{equation}
H^{-1}=1-\sum_{n\ge 0}(-1)^nR_{12^n}\,,
\end{equation}
which can again be interpreted as an identity in $\FQSym$
\begin{equation}
H^{-1}=1+\sum_{n\ge 0}(-1)^n\sum_{C(\sigma^{-1})=(12^n)}\F_\sigma\,.
\end{equation}


\end{document}